\documentclass[12pt]{article}
\usepackage{amsfonts,amssymb,amsmath}
\usepackage{graphicx}
\usepackage{color}
\usepackage{tikz}
\usepackage{lineno,hyperref}
\usepackage{comment}
\usepackage{float}

\title{\bf Prime Square Order Cayley Graph of Cyclic Groups of Particular Valency}
\author{Iqbal Adisna Atmaja$^{1,a}$,  Ahmad Erfanian$^{2,b}$,   Yeni Susanti$^{1,c}$,\\ Muhammad Nurul Huda$^{1,d}$, Ari Suparwanto$^{1,e,*}$ }

\date{}
		
\newcommand{\qed}{{\unskip\nobreak\hfil\penalty50\hskip .001pt \hbox{}
          \nobreak\hfil
          \vrule height 1.2ex width 1.1ex depth -.1ex
           \parfillskip=0pt\finalhyphendemerits=0\medbreak} }

\newenvironment{Proof}{\begin{trivlist}\item[\hskip%
\labelsep{\bf Proof.\quad}]}%
{\hfill\qed\end{trivlist}}
\newtheorem{Theorem}{Theorem}[section]
\newtheorem{Proposition}[Theorem]{Proposition}
\newtheorem{Corollary}[Theorem]{Corollary}
\newtheorem{Lemma}[Theorem]{Lemma}
       \newtheorem{rem}[Theorem]{Remark}
\newenvironment{Remark}{\begin{rem}\rm}{\end{rem}}
\newtheorem{define}[Theorem]{Definition} 

      \newtheorem{exer}[Theorem]{Exercise}

      \newtheorem{ex}[Theorem]{Example}
\newenvironment{Example}{\begin{ex}\rm}{\end{ex}}

\begin{document}
\maketitle
\begin{center}
{\small$^{1}$Department of Mathematics, Faculty of Mathematics and Natural Sciences, Universitas Gadjah Mada, Yogyakarta, Indonesia\\
$^{2}$Department of Pure Mathematics and Center of Excellence in Analysis on Algebraic Structures, Ferdowsi University of Mashhad, Mashhad, Iran}
\\

{\small $^{a}$ iqbaladisnaatmaja@mail.ugm.ac.id}\\
{\small $^{b}$} erfanian@um.ac.ir\\
{\small $^{c}$ yeni\_math@ugm.ac.id}\\
{\small $^{d}$ muhammadnurulhuda2000@mail.ugm.ac.id}\\
{\small $^{e}$ ari\_suparwanto@ugm.ac.id (corresponding author)}
\end{center}
\begin{abstract}
As a vital link between group theory and graph theory, Cayley graphs provide a geometric framework for encoding algebraic structures. This study explores the properties of Cayley graphs derived from cyclic groups whose order is the square of the product of three distinct prime numbers. We specifically examine cases where the connecting set is defined by the collection of all elements with an order equal to the square of a prime. A comprehensive analysis of these graphs is presented, focusing on structural characteristics such as connectivity, Eulerian properties, and Hamiltonicity. Furthermore, we determine several key graph parameters, including the clique number, chromatic number, independence number, and diameter.

\vspace{0.5cm}
\end{abstract}

\section{Introduction}

Arthur Cayley (1878) in \cite{cayley} first introduced the concept of the Cayley graph as a way to represent a group through a graph structure. The idea was later revisited under the term *“Gruppenbild”* (group diagram) in Max Dehn’s unpublished lectures on group theory. This development laid the groundwork for modern geometric group theory. Dehn employed a generating set to construct a geometric representation of groups, allowing them to be studied from a geometric perspective. This viewpoint led to the emergence of highly symmetric graphs, known as transitive graphs, which play an important role in both graph theory and group theory.

Cayley graphs are fundamental in tackling a range of problems, including Hamiltonian paths and cycles, modeling interconnection networks, studying group width, and constructing expander graphs—many of which have significant applications in computer science and related areas. Moreover, they serve as important instruments in both combinatorial and geometric group theory, where their properties are investigated from the perspectives of graph theory as well as algebra.

Let $G$ be any group with identity element $e_G$. Let $C$ be a nonempty subset of $G$ satisfying $C^{-1} \subseteq C$, meaning $C$ is closed under inversion and $e_G \notin C$. The Cayley graph of $G$ with respect to the connecting set $C$, denoted $Cay(G, C)$, is defined as a simple undirected graph with vertex set $G$, where two different vertices $x$ and $y$ are adjacent precisely when $xy^{-1} \in C$.
It is well-known that the Cayley graph is $|C|$-regular and vertex transitive. Moreover, the necessary and sufficient for the graph to be connected is that $C$ is the generator of $G$. 

While most research on Cayley graphs focuses on fixed finite groups, the connecting sets often vary. For example, some studies explore different choices of generating sets for specific groups. Additionally, some combinatorialists have worked to promote the study of Cayley graphs in more random cases, where the generating sets are chosen randomly, leading to new insights and results. B. Tolue \cite{tolue1, tolue2} defined prime order Cayley graphs, which are Cayley graphs of a group $G$ with set $C$ consisting of all non-identity elements of prime order in $G$. She also explored certain properties of these graphs, including the planarity of the prime order Cayley graph of abelian groups. Subsequently, I. Shojaee et al. \cite{shojaee} investigated and derived several graph-theoretic properties of Cayley graphs of prime order associated with groups of the form $\mathbb{Z}_{p^{n_1}} \times \cdots \times \mathbb{Z}_{p^{n_t}}$.
In \cite{YS}, the authors investigated Cayley graphs where the subset $C$ consists of prime square order elements instead of prime order elements. These are referred to as prime square order Cayley graphs.  The study in \cite{YS} focused on cyclic groups of order $p^2 q^2$ for distinct primes $p$ and $q$, a scope that was later expanded to cubic order elements in \cite{supar}. Both works analyzed various graph invariants, including the independence, chromatic, domination, and clique numbers.

In this paper, we investigate the prime square order Cayley graph for cyclic groups of order $\alpha^2 \beta^2 \gamma^2$, where $\alpha, \beta$, and $\gamma$ are distinct primes such that $\alpha < \beta < \gamma$. While our objectives align with the characterizations presented in \cite{YS} and \cite{supar}, this work employs a different analytical approach. Specifically, we utilize the cyclic decomposition of the group to derive the graph's properties, providing a more structural perspective on its characterization.

We begin with a brief overview of several notions from graph theory, followed by a theorem concerning Cayley graphs. Standard definitions in graph theory and group theory can be found in \cite{Graph theoretical concepts 2, Group theoretical concepts, west}. Let $\Gamma$ be a finite graph with vertex set $V(\Gamma)$ and edge set $E(\Gamma)$. A graph $H$ is called a subgraph of $\Gamma$ if $V(H) \subseteq V(\Gamma)$ and $E(H) \subseteq E(\Gamma)$. Furthermore, $H$ is said to be an induced subgraph if, for any vertices $x, y \in V(H)$, the presence of an edge $xy$ in $E(\Gamma)$ guarantees that $xy$ also belongs to $E(H)$; in this situation, $H$ is induced by $V(H)$. The degree of a vertex $v$ is defined as the number of edges in $E(\Gamma)$ that are incident to $v$. A graph $\Gamma$ is termed regular if all its vertices share the same degree; in particular, it is called $r$-regular when each vertex has degree $r$.

A trail in a graph $\Gamma$ is defined as a sequence of vertices and distinct edges of the form
$u_1 - u_1u_2 - u_2 - \dots - u_{k-1}u_k - u_k$, where $u_i \in V(\Gamma)$ for each $i = 1, 2, \dots, k$ and every edge $u_i u_{i+1}$ belongs to $E(\Gamma)$. The trail is said to be closed if its initial and terminal vertices are the same. A path is a trail in which all vertices are distinct, except possibly the first and last. A cycle is a path that is closed. An Eulerian trail is a closed trail that traverses each edge exactly once, and a graph admitting such a trail is called an Eulerian graph. A graph is termed semi-Hamiltonian if it contains a path that visits every vertex exactly once. Meanwhile, a Hamiltonian graph is one that includes a Hamiltonian cycle, that is, a cycle passing through every vertex exactly once.

The length of a path is defined as the number of edges it includes. A graph $\Gamma$ is said to be connected if for every pair of vertices, there exists a path joining them. In a connected graph, the distance between two vertices $x$ and $y$ is the length of a shortest path connecting them. The diameter of $\Gamma$ is the maximum distance among all pairs of distinct vertices.

A graph is called $k$-partite if its vertex set can be partitioned into $k$ pairwise disjoint nonempty subsets such that every edge joins vertices from different subsets. In the special case $k=2$, the graph is referred to as bipartite. It is said to be complete bipartite if every vertex in one subset is adjacent to all vertices in the other subset. An independent set in $\Gamma$ is a nonempty subset $I \subseteq V(\Gamma)$ whose vertices are pairwise nonadjacent. The maximum cardinality of such a set is called the independence number, denoted by $\alpha(\Gamma)$. A clique is a subset of vertices that induces a complete subgraph, and its maximum size is known as the clique number, denoted by $\omega(\Gamma)$. The chromatic number $\chi(\Gamma)$ is defined as the smallest number of colors needed to assign to the vertices so that adjacent vertices receive different colors. A dominating set $D \subseteq V(\Gamma)$ is a set such that every vertex in $V(\Gamma)\setminus D$ is adjacent to at least one vertex in $D$. The minimum size of such a set is called the domination number, denoted by $\gamma(\Gamma)$. A graph is planar if it can be drawn in the plane without any edge crossings, except at common endpoints. An automorphism of a graph $\Gamma$ is a bijection $f : V(\Gamma) \to V(\Gamma)$ that preserves adjacency in two directions, that is, $xy \in E(\Gamma)$ if and only if $f(x)f(y) \in E(\Gamma)$ for all vertices $x, y$. The graph is said to be vertex transitive if, for any two vertices $a$ and $b$, there exists an automorphism $f$ such that $f(a) = b$.

Consider a group $G$ and the set $C=\{x\in G| |x|=\mu^2, \mu ~{\rm is ~prime}\}$ where $|x|$ represents the order of $x$. Clearly, the identity element $e_G$ of $G$ does not belong to $C$  and $C^{-1}\subseteq C$ as $x$ and $x^{-1}$ are of the same order. 
The graph Cayley of $G$ respect to the connecting set $C$  will be denoted as $Cay_{p^2} (G,C)$. Obviously, $Cay_{p^2}(G,C)$ is the Cayley graph with $G$ as the vertex set  and two vertices $x$ and $y$ in $G$ are connected by an edge exactly when $xy^{-1}\in C$. Specifically, for a cyclic group $G$  of order $\alpha^2\beta^2\gamma^2$ for some distinct prime numbers $\alpha$, $\beta$ and $\gamma$, we have
$C=\{x\in G| |x|=\alpha^2 ~{\rm or}~ \beta^2~{\rm or}~\gamma^2\}$.  Throughout this manuscript, the group $G$ is cyclic, generated by $a\in G$ and of order $\alpha^2\beta^2\gamma^2$, where $\alpha$, $\beta$, and  $\gamma$ are distinct prime numbers with $\alpha<\beta<\gamma$. Moreover,  $C=\{x\in G||x|=\alpha^2 ~{\rm  or}~\beta^2~{\rm or}~\gamma^2\}$.

\section{Results}

Suppose that we have a group $G$ and the set $C=\{x\in G| |x|=\mu^2, \mu ~{\rm is ~prime}\}$ where $|x|$ represents the order of $x$. Clearly,    the identity element $e_G$ of $G$ does not belong to $C$  and $C^{-1}\subseteq C$ as $x$ and $x^{-1}$ are of the same order. The graph Cayley of $G$ respect to the connecting set $C$  will be denoted as $Cay_{p^2} (G,C)$. Obviously, $Cay_{p^2}(G,C)$ is the Cayley graph with  vertex set  $G$ and  for every $x$ and $y$ in $V(G)$, $xy\in E(G)$ exactly when $xy^{-1}\in C$. Specifically, for a cyclic group $G$  of order $\alpha^2\beta^2\gamma^2$ for some distinct prime numbers $\alpha$, $\beta$ and $\gamma$, we have
$C=\{x\in G| |x|=\alpha^2 ~{\rm or}~ \beta^2~{\rm or}~\gamma^2\}$.  Throughout this manuscript, the group $G$ is cyclic, generated by $a\in G$ and of order $\alpha^2\beta^2\gamma^2$, where $\alpha$, $\gamma$, and  $\beta$ are distinct prime numbers with $\alpha<\beta<\gamma$. Moreover,  $C=\{x\in G||x|=\alpha^2 ~{\rm  or}~\beta^2~{\rm or}~\gamma^2\}$.

We have the following preliminary results as attempted in \cite{YS}.

\begin{Lemma}\label{Lem-C}
These two assertions hold.

\begin{enumerate}
\item[\mbox{\rm (i)}] $C=\{a^{k\beta^2\gamma^2}|k\in\{1,\ldots,\alpha^2-1\}\setminus\{\alpha l|1\leq l\leq \alpha-1\}\}\cup \{a^{k\alpha^2\gamma^2}|k\in\{1,\ldots,\beta^2-1\}\setminus\{\beta l|1\leq l\leq \beta-1\}\}\{a^{k\alpha^2\beta^2}|k\in\{1,\ldots,\gamma^2-1\}\setminus\{\gamma l|1\leq l\leq \gamma-1\}\}$
\item[\mbox{\rm (ii)}] $|C|=\alpha^2+\beta^2+\gamma^2-\alpha-\beta-\gamma$.
\end{enumerate}
\end{Lemma}
\begin{Proof}
\begin{enumerate}
\item[\mbox{\rm (i)}] Obviously, $|G|=|<a>|=\alpha^2\beta^2\gamma^2$. The set $C$ includes all elements in $G$ having order $\alpha^2$ or $\beta^2$ or $\gamma^2$. The elements of order $\alpha^2$ are of the form
$a^{k\beta^2\gamma^2}$, where $1\leq k\leq \alpha^2-1$ but $k\neq \alpha l$, for $l=1,2,\ldots,\alpha-1$, the elements of order $\beta^2$ are of the form
$a^{k\alpha^2\gamma^2}$ where $1\leq k\leq \beta^2-1$ and $k\neq \beta l$, for $l=1,2,\ldots,\beta-1$, and the elements of order $\gamma^2$ are of the form $a^{k\alpha^2\beta^2}$ where $1\leq k\leq \gamma^2-1$ and $k\neq \gamma l$, for $l=1,2,\ldots,\gamma-1$. Hence, $C=\{a^{k\beta^2\gamma^2}|k\in\{1,\ldots,\alpha^2-1\}\setminus\{\alpha l|1\leq l\leq \alpha-1\}\}\cup \{a^{k\alpha^2\gamma^2}|k\in\{1,\ldots,\beta^2-1\}\setminus\{\beta l|1\leq l\leq \beta-1\}\}\{a^{k\alpha^2\beta^2}|k\in\{1,\ldots,\gamma^2-1\}\setminus\{\gamma l|1\leq l\leq \gamma-1\}\}$. 
\item[\mbox{\rm (ii)}] From (i) we know that $|C|=\alpha^2+\beta^2+\gamma^2-\alpha-\beta-\gamma$.
\end{enumerate} 
\end{Proof}

Now, we will write any element $x$ of $G$ in the form $x=a^{i+j\alpha^2 +k\alpha^2\beta^2}$ where $0\leq i\leq \alpha^2-1$, $0\leq j\leq \beta^2-1$ and $0\leq k\leq \gamma^2-1$. Clearly, 
$G=\{a^{i+j\alpha^j+k\alpha^2\beta^2}|0\leq i\leq \alpha^2-1, 0\leq j\leq \beta^2-1, 0\leq k\leq \gamma^2-1\}$. Let also classify the elements of $G$ into the following sets. For a fixed $r$, $0\leq r\leq \alpha^2-1$, a fixed $s$, $0\leq s\leq \beta^2-1$, and a fixed $t$, $0\leq t\leq \gamma^2-1$, let
$$R_{r}=\{a^{r+j\alpha^2+k\alpha^2\beta^2}|0\leq j\leq \beta^2-1~{\rm and}~ 0\leq k\leq \gamma^2-1\}$$
$$S_{s}=\{a^{i+s\alpha^2+k\alpha^2\beta^2}|0\leq i\leq \alpha^2-1 ~{\rm and}~ 0\leq k\leq \gamma^2-1\}$$
$$T_{t}=\{a^{i+j\alpha^2+t\alpha^2\beta^2}|0\leq i\leq \alpha^2-1 ~{\rm and}~ 0\leq j\leq \beta^2-1\}.$$ Clearly, any element of the form $a^{i+j\alpha,.m2+k\alpha^2\beta^2}$ is in $R_i\cap S_j\cap T_k$. Moreover, we obtain the following properties.

\begin{Proposition} \label{Prop2.3}The statements below are valid.k
\begin{enumerate}
    \item[\mbox{\rm (i)}] For any $t$, $0\leq t\leq \gamma^2-1$, it follows that $T_t$ is an  independent set. 
    \item[\mbox{\rm (ii)}] For any $r$ and $s$, where $0\leq r\leq \beta^2-1$ and $0\leq s\leq \alpha^2-1$, and for any $k$ and $k'$ where $0\leq k,k'\leq \gamma^2-1$, it follows that two different elements $x=a^{r+s\alpha^2+k\alpha^2\beta^2}$ and $y=a^{r+s\alpha^2+k'\alpha^2\beta^2}$ in $R_r\cap S_s$ are adjacent if and only if $k\not\equiv k' (mod~\gamma)$.   
    \item[\mbox{\rm (iii)}] For any $r$ and $s$ with $0\leq r\leq \alpha^2-1$, $0\leq s\leq \beta^2-1$, the  subgraph induced by $R_r\cap S_s$ contains a cycle.
    \item[\mbox{\rm (iv)}] For any $r$, $0\leq r\leq \alpha^2-1$, and any $s$, $0\leq s\leq \beta^2-1$ but $(r,s)\neq (0,0)$, it follows that $|(R_r\cap S_s)\cap \{a^{k\gamma^2}|k=1,2,\ldots,\alpha^2\beta^2-1\}|=1$.

\item[\mbox{\rm (v)}] Let for any $r=0,1,\ldots, \alpha^2-1$, $C_{\beta^2}^r=\{x\in G|x=a^{k\alpha^2\gamma^2+r\gamma^2}, 0\leq k\leq \beta^2-1\}$. Then, we have $C_{\beta^2}^r\subseteq R_r$.

\item[\mbox{\rm (vi)}] Let $C_{\alpha^2}=\{x\in G|x=a^{k\beta^2\gamma^2}, 0\leq k\leq \alpha^2-1\}$. For any $r=0,1,\ldots,\alpha^2-1$, it follows that $|R_r\cap C_{\alpha^2}|=1$.

\item[\mbox{\rm (vii)}] For any $r$, $0\leq r\leq \alpha^2-1$, there exists $x_r\in R_r$ such that $\{x_r|0\leq r\leq \alpha^2-1\}$ forms a cycle.

\item[\mbox{\rm (viii)}] For any $r$, $0\leq r\leq \alpha^2-1$ and for any $s$, $0\leq s\leq \beta^2-1$ there exists $x_{r,s}\in R_r\cap S_s$ such that $\{x_{r,s}|0\leq s\leq \beta^2-1\}$ forms a cycle.

\end{enumerate}
    \end{Proposition}
\begin{Proof}
    \begin{enumerate}
    \item[\mbox{\rm (i)}] Let $t$ with $0\leq t\leq \gamma^2-1$  and let $x,y\in T_t$ be arbitrary. Let $x=a^{i+j\alpha +t\alpha^2\beta^2}$ and $y=a^{i'+j'\alpha+t\alpha^2\beta^2}$ such that $(i,j)\neq (i',j')$. It follows that $xy^{-1}=a^{(i-i')+(j-j')\alpha}\notin C$. Thus, $x$ and $y$ are not connected by any single edge. As a consequence, $T_t$ is an independence set. 

    \item[\mbox{\rm (ii)}] Let $x=a^{r+s\alpha^2+k\alpha^2\beta^2}$ and $y=a^{r+s\alpha^2+k'\alpha^2\beta^2}$ where  $0\leq r\leq \alpha^2-1$, $0\leq s\leq \beta^2-1$
    and $0\leq k,k'\leq \gamma^2-1$. Then we have $xy^{-1}=a^{(k-k')\alpha^2\beta^2}.$ If $x$ and $y$ are adjacent, then $xy^{-1}\in C$ and hence $k-k'$ is not in the form $l\gamma$ for some $1\leq l\leq \gamma-1$. This means that $k\not\equiv k'~(mod ~\gamma)$. Conversely, if $k\not\equiv k'~(mod ~\gamma)$, then $k-k'\neq l\gamma$ for some $l$. Thus, $x$ and $y$ are not connected by any single edge.

    \item[\mbox{\rm (iii)}] By (ii), we have a cycle $a^{r+s\alpha^2}-a^{r+s\alpha^2+\alpha^2\beta^2}-a^{r+s\alpha^2+2\alpha^2\beta^2}-\ldots-a^{r+s\alpha^2+(\gamma^2-1)\alpha^2\beta^2}-a^{r+s\alpha^2}$ in the subgraph induced by $R_r\cap S_s$.

    \item[\mbox{\rm (iv)}] Let $A=\{a^{k\gamma^2}|k=1,2,\ldots,\alpha^2\beta^2-1\}$. Assume that there are some $r$ and $s$, $0\leq r\leq \alpha^2-1$, $s$, $0\leq s\leq \beta^2-1$ and $(r,s)\neq (0,0)$, such that $R_r\cap S_s$ contains two elements in $A$. Let the two elements are $u=a^{k\gamma^2}$ and $v=a^{l\gamma^2}$. WLOG, let $l<k$. As $u,v\in R_r\cap S_s$, let also
$u=a^{r+s\alpha^2+p\alpha^2\beta^2}$
and $v=a^{r+s\alpha^2+q\alpha^2\beta^2}$ with $p\neq q$.  Then we have $uv^{-1}=a^{k-l}\gamma^2\in R_0\cap S_0$ and $uv^{-1}=a^{(p-q)\alpha^2\beta^2}\not\in R_0\cap S_0$. Thus we have a contradiction. Therefore, for each $(r,s)\neq (0,0)$, the set $R_r\cap S_s$ contains at most one element from $A$. But, the absence of elements $A$ in $R_r\cap S_s$ implies that there is $(r,s)\neq (0,0)$ with $R_r\cap S_s$ contains two elements from $A$. Therefore, we conclude that $R_r\cap S_s\cap A$ is a singleton. This completes the proof.

\item[\mbox{\rm (v)}] Let $x\in C_{\beta^2}^r$, i.e. $x=a^{k\alpha^2\gamma^2+r\gamma^2}$ for some $k$ with $0\leq k\leq \beta^2-1$. Assume that $x\in R_i$, $i\neq r$. Then, simultaneously we have $x=a^{i+j\alpha^2+k'\alpha^2\beta^2}$, $i\neq r$,  $0\leq j\leq  \beta^2-1$ and $0\leq k'\leq \gamma^2-1$. We have $$a^{i-r\gamma^2+\alpha^2(j+k'\beta^2-k\gamma^2)}=e_G.$$As a consequence, $\alpha^2\beta^2\gamma^2$ divides $i-r\gamma^2+\alpha^2(j+k'\beta^2-k\gamma^2)$. This can happen only if $i=0$ and $r=0$. But then we have a contradiction as $i\neq r$. Therefore, $C_{\beta^2}^r\subseteq R_r$.

\item[\mbox{\rm (vi)}] Clearly, $a^0\in R_0\cap C_{\alpha^2}$. 
Let assume that there exists $r$, $0\leq r\leq \alpha^2-1$, such that $|R_r\cap C_{\alpha^2}|\geq 2$. Let $u,v\in R_r\cap C_{\alpha^2}$, be two different elements. Let $u=a^{r+j\alpha^2+k\alpha^2\beta^2}$ and $v=a^{r+j'\alpha^2+k'\alpha^2\beta^2}$ and let assume simultaneously $u=a^{l\beta^2\gamma^2}$ and $v=a^{l'\beta^2\gamma^2}$, with $0\leq j,j'\leq \beta^2-1$, $0\leq k,k'\leq \gamma^2-1$, $0\leq l,l'\leq \alpha^2-1$.

We obtain 
$uv^{-1}=a^{(j-j')\alpha^2+(k-k')\alpha^2\beta^2}=a^{(l-l')\beta^2\gamma^2}$.
We have $$a^{(j-j')\alpha^2+((k-k')\alpha^2-(l-l')\gamma^2)\beta^2}=e_G.$$
As $\alpha^2\beta^2\gamma^2$ does not divide  ${(j-j')\alpha^2+((k-k')\alpha^2-(l-l')\gamma^2)\beta^2}$, necessarily ${(j-j')\alpha^2+((k-k')\alpha^2-(l-l')\gamma^2)\beta^2}=0$. Thus,  $j=j', k=k'$, and $l=l'$.

\item[\mbox{\rm (vii)}] By (vi), for each $r$, $1\leq r\leq \alpha^2-1$, there exists precisely one $x_r\in R_r\cap C_{\alpha^2}$. This $x_r$'s are of the form $a^{l\beta^2\gamma^2}$, $1\leq l\leq \alpha^2-1$. Now, let $l_1<\ldots\leq l_{\alpha^2-1}$.  As $0\in R_0$, we then have a cycle of the form $0-a^{l_1\beta^2\gamma^2}-a^{l_2\beta^2\gamma^2}-\ldots-a^{l_{\alpha^2-1}\beta^2\gamma^2}-0$. This completes the assertion (vii).

\item[\mbox{\rm (viii)}] 
By (vi), for any $r=1,2,\ldots,\alpha^2-1$, there exists a unique element in $R_r\cap C_{\alpha^2}$ of the form $a^{k_r\beta^2\gamma^2}$.  Clearly, for any $l=0,1,\ldots,\beta^2-1$, $a^{k_r\beta^2\gamma^2+l\alpha^2\gamma^2}\in R_r$. Thus, $A_r=\{a^{k_r\beta^2\gamma^2+l\alpha^2\gamma^2}|0\leq l\leq \beta^2-1\}\subseteq R_r$. Moreover, from the set $A_r$, we can form a cycle $a^{k_r\beta^2\gamma^2} -  a^{k_r\beta^2\gamma^2+\alpha^2\gamma^2}-\cdots-a^{k_r\beta^2\gamma^2+(\beta^2-1)\alpha^2\gamma^2}$. Now, we will show that $|A_r\cap S_s|=1$ for each $s$, where $0\leq s\leq \beta^2-1$. Assume that there are two distinct elements $a^{k_r\beta^2\gamma^2+l\alpha^2\gamma^2}$ and $a^{k_r\beta^2\gamma^2+l'\alpha^2\gamma^2}$ in $A_r\cap S_s$. We have  $a^{k_r\beta^2\gamma^2+l\alpha^2\gamma^2}=a^{r+s\alpha^2+j\alpha^2\beta^2}$ and $a^{k_r\beta^2\gamma^2+l'\alpha^2\gamma^2}=a^{r+s\alpha^2+j'\alpha^2\beta^2}$ for some $0\leq j, j'\leq \gamma^2-1$. As a consequence, we have
$a^{k_r\beta^2\gamma^2+l\alpha^2\gamma^2-r-s\alpha^2-j\alpha^2\beta^2}=e_G$ and $a^{k_r\beta^2\gamma^2+l'\alpha^2\gamma^2-r-s\alpha^2-j'\alpha^2\beta^2}=e_G$ so that $k_r\beta^2\gamma^2+l\alpha^2\gamma^2-r-s\alpha^2-j\alpha^2\beta^2=m\alpha^2\beta^2\gamma^2$ and $k_r\beta^2\gamma^2+l'\alpha^2\gamma^2-r-s\alpha^2-j'\alpha^2\beta^2=m'\alpha^2\beta^2\gamma^2$ for some integers $m$ and $m'$. Hence, we have $(l-l')\alpha^2\gamma^2-(j-j')=(m-m')\alpha^2\beta^2\gamma^2$.
It can happen only when $l=l'$, $j=j'$, and $m=m'$, meaning that $a^{k_r\beta^2\gamma^2+l\alpha^2\gamma^2}=a^{k_r\beta^2\gamma^2+l'\alpha^2\gamma^2}$. Thus, we have a contradiction. We can conclude that $|A_r\cap S_s|=1$. This proves assertion (viii).

\end{enumerate}\end{Proof}




\begin{Proposition}
The graph $Cay_{p^2}(G,C)$ is not planar and has girth $3$.
\end{Proposition}
\begin{Proof}
   By Proposition \ref{Prop2.3} part (ii), all vertices of the forms $a^{r+s\alpha^2+i\alpha^2\beta^2}$, $i=0,1,\ldots,\gamma-1$,  induce the complete subgraph $K_{\gamma}$ in $Cay_{p^2}(G,C)$. As $\alpha<\beta<\gamma$ which implies  $\gamma\geq 5$,   $K_{\gamma}$ contains the complete subgraph $K_{5}$. By Kuratowski's theorem, we conclude that $Cay_{p^2}(G,C)$ is not planar and hence the girth of $Cay_{p^2}(G,C)$ is $3$.
    \end{Proof}

Now, we will explore further properties of the Cayley graph using cyclic group decomposition. Let $G = \langle a \rangle $ be a cyclic group of order $\alpha^2 \beta^2 \gamma^2$, where $\alpha < \beta < \gamma$ are distinct primes. Then $G \cong \mathbb{Z}_{\alpha^2 \beta^2 \gamma^2}$. By the cyclic group decomposition, we further have
    \begin{align*}
        \mathbb{Z}_{\alpha^2 \beta^2 \gamma^2} \cong \mathbb{Z}_{\alpha^2} \times \mathbb{Z}_{\beta^2} \times \mathbb{Z}_{\gamma^2}.
    \end{align*}
An explicit isomorphism $f : G \to \mathbb{Z}_{\alpha^2} \times \mathbb{Z}_{\beta^2} \times \mathbb{Z}_{\gamma^2}$ is given by
    \begin{align*}
        f(a^k) = (k \text{ (mod }\alpha^2), k \text{ (mod }\beta^2), k \text{ (mod }\gamma^2)),
    \end{align*}
for all $g = a^k \in G$.

\begin{Remark}
    Let $g=a^k \in G$. We define $g_\alpha := k \text{ (mod }\alpha^2)$, $g_\beta := k \text{ (mod }\beta^2)$, and $g_\gamma := k \text{ (mod }\gamma^2)$, which are referred to as the $\alpha$-component, $\beta$-component, and $\gamma$-component of $g$, respectively. In other words, every element $g \in G$ can be represented as the ordered triple $(g_\alpha, g_\beta, g_\gamma)$.
\end{Remark}

\begin{Proposition}\label{thconnected}
    The Cayley graph $Cay_{p^2}(G,C)$ is connected. 
\end{Proposition}

\begin{Proof}
    Observe that $a^{\beta^2 \gamma^2}, a^{\alpha^2 \gamma^2}, a^{\alpha^2\beta^2 } \in C$. Since $\gcd (\alpha^2 \gamma^2, \alpha^2 \gamma^2, \alpha^2\beta^2) = 1$, Bézout's Identity guarantees that there exist integers $u,v$, and $w$ such that
        \begin{align*}
            u(\alpha^2 \gamma^2)+ v(\alpha^2 \gamma^2) + w( \alpha^2\beta^2) = 1.
        \end{align*}
    Consequently, 
        \begin{align*}
            a^1 = a^{u(\alpha^2 \gamma^2)+ v(\alpha^2 \gamma^2) + w( \alpha^2\beta^2)} = a^{u(\alpha^2 \gamma^2)}a^{v(\alpha^2 \gamma^2)}a^{w( \alpha^2\beta^2)} \in \langle C \rangle 
        \end{align*}
    since each factor belongs to $\langle C \rangle$. Hence $G = \langle a \rangle \subseteq \langle C \rangle$, and therefore $C$ generates $G$. It follows that the Cayley graph $Cay_{p^2}(G,C)$ is connected.
\end{Proof}

As a consequence of Theorem \ref{thconnected} and Lemma \ref{Lem-C}, we have the following property on the regularity and eulerianity of $Cay_{p^2}(G,C)$.

\begin{Proposition}
The graph $Cay_{p^2}(G,C)$ is $(\alpha^2+\beta^2+\gamma^2-\alpha-\beta-\gamma)$-regular and eulerian.
\end{Proposition}
\begin{Proof}
    By the fact that the Cayley graph is $|C|$-regular and by Lemma \ref{Lem-C}, it is obvious that $Cay_{p^2}(G,C)$ is $(\alpha^2+\beta^2+\gamma^2-\alpha-\beta-\gamma)$-regular. Moreover, as $\alpha^2+\beta^2+\gamma^2-\alpha-\beta-\gamma$  is an even number, we conclude that $Cay_{p^2}(G,C)$ is eulerian. 
\end{Proof}
Beyond the global regularity of the graph, it is essential to investigate its local density, specifically the existence of complete subgraphs. By analyzing the adjacency conditions between the vertices within the cyclic decomposition $\mathbb{Z}_{\alpha^2} \times \mathbb{Z}_{\beta^2} \times \mathbb{Z}_{\gamma^2}$, we can determine the maximum number of vertices that are mutually adjacent. The following proposition establishes the clique number of the graph, which is fundamentally bounded by the largest prime factor of the group's order.

\begin{Proposition}
    The clique number of the Cayley graph $Cay_{p^2}(G,C)$ is $\gamma$.
\end{Proposition}

\begin{Proof}
    Let $g=(g_\alpha, g_\beta, g_\gamma)$ and $h=(h_\alpha, h_\beta, h_\gamma)$ be adjacent vertices in $V(Cay_{p^2}(G,C))$. By definition of adjacency, the order of $gh^{-1}$ is equal to $\alpha^2$, $\beta^2$, or $\gamma^2$. If $|gh^{-1}|= \alpha^2$, then $g_\alpha \not\equiv h_\alpha \text{ (mod }\alpha)$, while $g$ and $h$ coincide in their $\beta-$ and $\gamma-$components. The cases $|gh^{-1}|=\beta^2$ or $\gamma^2$ follow by analogous arguments. Consequently, any clique set $T$ must lie entirely in one of the following sets, 
        \begin{align*}
            P = \{(x_1, y_1, z_1) : x_1& \in \mathbb{Z}_{\alpha^2} \}, \hspace{0.2cm} Q = \{(x_2, y_2, z_2) : y_2 \in \mathbb{Z}_{\beta^2} \}, \text{ or } \\
            R &= \{(x_3, y_3, z_3) : z_3 \in \mathbb{Z}_{\gamma^2} \}.
        \end{align*}
    Suppose that $T$ is contained in $P$. Identifying $P$ as a subset of $\mathbb{Z}_{\alpha^2}$, two vertices in $T$ are adjacent if and only if they belong to distinct residue class modulo $\alpha$. Hence, the largest possible clique in $P$ is
        \begin{align*}
            T = \{ \bar{0}, \bar{1}, \dots , \overline{\alpha-1} \} \subseteq \mathbb{Z}_{\alpha^2},
        \end{align*}
    which has size $\alpha$. Similarly, if $T \subseteq Q$ or $T \subseteq R$, the corresponding clique numbers are $\beta$ and $\gamma$, respectively. Therefore, the clique number of $Cay_{p^2}(G,C)$ is $\max \{\alpha, \beta, \gamma \}= \gamma$.
\end{Proof}

The determination of the clique number provides a fundamental lower bound for the chromatic number of the graph, as $\chi(Cay_{p^2}(G,C)) \geq \omega(Cay_{p^2}(G,C)) = \gamma$. To show that this bound is tight, we construct an explicit vertex coloring using the properties of the cyclic decomposition. By partitioning the vertex set based on the residue classes modulo the prime factors, we can ensure that no two adjacent vertices share the same color. This leads to the following result regarding the chromatic number.

\begin{Proposition}
    The chromatic number of the Cayley graph $Cay_{p^2}(G,C)$ is $\gamma$.
\end{Proposition}

\begin{Proof}
    Since the clique number of $Cay_{p^2}(G,C)$ is $\gamma$, it follows that $\chi(Cay_{p^2}(G,C)) \geq \gamma$. Now, we wish to construct a proper coloring of $Cay_{p^2}(G,C)$ using exactly $\gamma$ colors. As $\alpha, \beta < \gamma$, there exist injective maps $\iota_\alpha : \mathbb{Z}_\alpha \hookrightarrow \mathbb{Z}_{\gamma}$ and $\iota_\beta: \mathbb{Z}_{\beta} \hookrightarrow \mathbb{Z}_{\gamma}$. For example, one may take
        \begin{align*}
            \iota_\alpha(k) = k \hspace{0.3cm}\text{  and  } \hspace{0.3cm} \iota_\beta(l) = l
        \end{align*}
    for $k =0,1, \dots, \alpha-1$ and $l =0,1,\dots, \beta -1$.

    Define a map $\zeta: V(Cay_{p^2}(G,C)) \to \mathbb{Z}_{\gamma}$ by 
        \begin{align*}
            \zeta((g_\alpha, g_\beta, g_\gamma)) \equiv \iota_\alpha( g_\alpha \text{ (mod }\alpha)) + \iota_\beta( g_\beta \text{ (mod }\beta)) + g_\gamma \hspace{0.3cm} \text{ (mod }\gamma).
        \end{align*}
    We show that $\zeta$ is a proper coloring. Let $g=(g_\alpha, g_\beta, g_\gamma)$ and $h=(h_\alpha, h_\beta, h_\gamma)$ be adjacent vertices. Suppose first that $|gh^{-1}| = \alpha^2$, the cases $|gh^{-1}| = \alpha^2$ or $\beta^2$ follow analogously. Then $g_\alpha \not\equiv h_\alpha \text{ (mod }\alpha)$, while $g_\beta = h_\beta$ and $g_\gamma = h_\gamma$. Hence, 
        \begin{align*}
            \zeta(g) - \zeta(h) \equiv \iota_\alpha(g_\alpha) - \iota_\alpha(h_\alpha) \hspace{0.3cm}\text{ (mod }\gamma).
        \end{align*}
    Since $\iota_\alpha$ is injective, we have $\iota_\alpha(g_\alpha) \not\equiv \iota_\alpha (h_\alpha) \text{ (mod }\gamma)$, and therefore $\zeta(g) \neq \zeta(h)$. Thus, $\zeta$ assigns distinct colors to adjacent vertices and is a proper $\gamma$-coloring of $Cay_{p^2}(G,C)$. Consequently, $\chi(Cay_{p^2}(G,C)) = \gamma$.
\end{Proof}

Having established the minimum number of colors required to partition the vertex set into independent sets, we now shift our focus to determining the maximum size of such a set, known as the independence number $\alpha(Cay_{p^2}(G,C))$. In the following propositions, we have some properties related to independence number and develop several technical tools to obtain the exact independence number of the graph.

\begin{Proposition}
    Let $(i,j,k) \in \mathbb{Z}_{\alpha} \times \mathbb{Z}_{\beta} \times \mathbb{Z}_{\gamma}$. The set $C(i,j,k) \subseteq V(Cay_{p^2}(G,C))$ defined as 
        \begin{align*}
            C(i,j,k) = \{(i+\alpha x, j+\beta y, k+ \gamma z): x \in \mathbb{Z}_{\alpha}, y \in \mathbb{Z}_{\beta}, z \in \mathbb{Z}_{\gamma} \}
        \end{align*}
    is an independent set of size $\alpha\beta\gamma$. 
\end{Proposition}

\begin{Proof}
    Let $g=(i+\alpha x, j+\beta y, k+ \gamma z)$ and $h= (i+\alpha x', j+\beta y', k+ \gamma z')$ be distinct vertices in $C(i,j,k)$. If $g$ and $h$ differ in two or three coordinates, then clearly they are non-adjacent. It remains to consider the case in which they differ in exactly one coordinate. Suppose that $x = x', y= y',$ and $z \neq z'$, the remaining cases follow symmetrically. Then $k+ \gamma z \equiv k+ \gamma z' \text{ (mod }\gamma)$, which implies that $|gh^{-1}| \neq \gamma^2$. Hence $g$ and $h$ are not adjacent. Therefore, any two distinct vertices of $C(i,j,k)$ are non-adjacent. 
\end{Proof}

\begin{Remark}
    The collection 
    \begin{align*}
        \mathcal{P}= \{C(i,j,k) : (i,j,k) \in \mathbb{Z}_{\alpha} \times \mathbb{Z}_{\beta} \times \mathbb{Z}_{\gamma} \}    
    \end{align*}
    forms a partition of the vertex set $V(Cay_{p^2}(G,C))$.
\end{Remark}

\begin{Proposition}
    Let $(i,j,k), (i',j',k') \in \mathbb{Z}_{\alpha} \times \mathbb{Z}_{\beta} \times \mathbb{Z}_{\gamma}$ be distinct indices. Then $C(i,j,k)$ and $C(i',j',k')$ are independent to one another if and only if no distinct indices agree on two coordinates.
\end{Proposition}

\begin{Proof}
    Let $(i,j,k), (i',j',k') \in \mathbb{Z}_{\alpha} \times \mathbb{Z}_{\beta} \times \mathbb{Z}_{\gamma}$ be distinct indices and let $g \in C(i,j,k)$ and $h \in C(i',j',k')$. We show that adjacency between $C(i,j,k)$ and $C(i',j',k')$ occurs if and only if the indices $(i,j,k)$ and $(i',j',k')$ agree on exactly two coordinates. 

    First, suppose that the indices agree in exactly two coordinates. Without the loss of generality, assume that $i \neq i'$, $j = j'$, and $k = k'$, the remaining cases follow by symmetry. Now, choose $u = (i+ \alpha x, j+ \beta y, k + \gamma z) \in C(i,j,k)$ and $v = (i'+ \alpha x, j'+ \beta y, k' + \gamma z) \in C(i',j',k')$. Then $uv^{-1} = ((i-i'), 0, 0)$. Since $i,i' \in \mathbb{Z}_{\alpha}$ with $i \not\equiv i' \text{ (mod }\alpha)$, it follows that $|uv^{-1}|= \alpha^2$. Hence, $u$ and $v$ are adjacent, and therefore adjacency occurs between $C(i,j,k)$ and $C(i',j',k')$.

    Conversely, suppose that the indices do not agree in two coordinates, we will show that no adjacency is possible. Let $g \in C(i,j,k)$ and $h \in C(i',j'.k')$. If the indices differ in all three coordinates, then the order of $gh^{-1}$ is divisible by $\alpha, \beta,$ and $\gamma$ simultaneously, and thus cannot be equal to $\alpha^2, \beta^2$, or $\gamma^2$. Hence, $g$ and $h$ are non-adjacent. 
    
    Next, assume that the indices agree in exactly one coordinate. Without loss of generality, let $i = i'$, $j \neq j'$, and $k \neq k'$. In this case, the order of $gh^{-1}$ is divisible by both $\beta^2$ and $\gamma^2$, and therefore cannot be equal to $\alpha^2, \beta^2$, or $\gamma^2$. Thus, no adjacency occurs. 
    
    Consequently, adjacency between $C(i,j,k)$ and $C(i',j',k')$ occurs if and only if the corresponding indices agree in exactly two coordinates. 
\end{Proof}

\begin{Corollary}\label{index restriction}
    Let $S$ be any independent set in $Cay_{p^2}(G,C)$ and let $\mathcal{I}$ be the set of indices defined as follows.
        \begin{align*}
            \mathcal{I} = \{ (i,j,k) \in \mathbb{Z}_{\alpha} \times \mathbb{Z}_{\beta} \times \mathbb{Z}_{\gamma} : S \cap C(i,j,k) \neq \varnothing \}.
        \end{align*}
    Then any pair of indices in $\mathcal{I}$ do not agree in two coordinates.
\end{Corollary}

\begin{Lemma}\label{maximum index}
    Let $\mathcal{I}$ be the set of indices defined in Corollary \ref{index restriction}. Then $\mathcal{I}$ has at most $\alpha \beta$ elements.
\end{Lemma}

\begin{Proof}
    Define the projection map $\pi_{\alpha \beta} : \mathcal{I} \to \mathbb{Z}_{\alpha} \times \mathbb{Z}_{\beta}$ as follows. 
        \begin{align*}
            \pi_{\alpha \beta} (i,j,k) = (i,j).
        \end{align*}
    We show that $\pi_{\alpha \beta}$ is an injective function. For the sake of contradiction, suppose that $\pi_{\alpha \beta}$ is not injective. This means there exist two distinct indices $a_1,a_2 \in \mathcal{I}$ such that $\pi_{\alpha \beta}(a_1) = \pi_{\alpha \beta} (a_2)$. Let $a_1 = (i,j,k)$ and $a_2 = (i',j',k')$. If $\pi_{\alpha \beta}(a_1) = \pi_{\alpha \beta} (a_2)$, then $i=i'$ and $j=j'$. This implies that $a_1$ and $a_2$ agree in at least two coordinates. By the definition of $\mathcal{I}$, no two elements can agree in exactly two coordinates. Therefore, $a_1$ and $a_2$ must also agree in the third coordinate, i.e., $k = k'$. This contradicts the assumption that $a_1$ and $a_2$ are distinct indices of $\mathcal{I}$. Thus, $\pi_{\alpha \beta}$ is injective.  

    Since $\pi_{\alpha \beta}$ is injective, the number of elements in $\mathcal{I}$ cannot exceed the number of possible outputs in the codomain $\mathbb{Z}_{\alpha} \times \mathbb{Z}_{\beta}$. Therefore, $|\mathcal{I}| \leq |\mathbb{Z}_{\alpha} \times \mathbb{Z}_{\beta}| = \alpha \beta$. By applying the same logic to projections $\pi_{\alpha \gamma}$ and $\pi_{\beta \gamma}$, defined similarly, we obtain $|\mathcal{I}| \leq \alpha \gamma$ and $|\mathcal{I}| \leq \beta \gamma$. Thus, $|\mathcal{I}| \leq \min \{\alpha \beta, \alpha \gamma, \beta \gamma \} = \alpha \beta$.  
\end{Proof}

\begin{Proposition}
    Let $S$ be any independent set in $Cay_{p^2}(G,C)$. Then $|S| \leq \alpha^2 \beta^2 \gamma$.
\end{Proposition}

\begin{Proof}
    Let $S$ be any independent set and let
    \begin{align*}
        \mathcal{I} = \{(i,j,k) \in \mathbb{Z}_{\alpha} \times \mathbb{Z}_{\beta} \times \mathbb{Z}_{\gamma} : S \cap C(i,j,k) \neq \varnothing \}.    
    \end{align*}
    Since $\mathcal{P}= \{C(i,j,k) : (i,j,k) \in \mathbb{Z}_{\alpha} \times \mathbb{Z}_{\beta} \times \mathbb{Z}_{\gamma} \}$ forms a partition on $V(Cay_{p^2}(G,C))$, we have
        \begin{align*}
            |S| = \sum_{(i,j,k) \in \mathcal{I}} |S \cap C(i,j,k)| \leq \sum_{(i,j,k) \in \mathcal{I}} \alpha \beta \gamma = |\mathcal{I}| \alpha \beta \gamma.
        \end{align*}
    Lemma \ref{maximum index} already establishes that $|\mathcal{I}| \leq \alpha \beta$, and therefore we obtain $|S| \leq \alpha^2 \beta^2 \gamma$. 
\end{Proof}

\begin{Lemma}\label{index in max construction}
    Define the set of indices $\mathcal{L}$ as follows
        \begin{align*}
            \mathcal{ L} = \{ (i,j,i+j \text{ (mod }\gamma)) : (i,j) \in \mathbb{Z}_{\alpha} \times \mathbb{Z}_{\gamma}  \}. 
        \end{align*}
    Then any pair of indices in $\mathcal{L}$ do not agree in exactly two coordinates. Moreover, it is clear that $|\mathcal{L}| = \alpha \beta$. 
\end{Lemma}

\begin{Proof}
    Let $a_1, a_2 \in \mathcal{L}$ be distinct indices. Let $a_1 = (i,j,k)$ and $a_2 = (i',j',k')$. Now, we examine that $a_1$ and $a_2$ do not share the same two coordinates. If $i = i'$ and $j = j'$, then it forces $k = k'$ making them the same element. If $i = i'$ and $k = k'$, then we have $i + j \equiv i' + j' \text{ (mod }\gamma)$, which implies that $j \equiv j' \text{ (mod }\gamma)$. Now, since $j, j' < \beta < \gamma$, the only way that they lie in the same residue class modulo $\gamma$ is if $j = j'$. This forces agreement in all three coordinates. Lastly, if $j=j'$ and $k =k'$, analogously, it also forces $i=i'$, making agreement in all three coordinates. In all cases where elements agree in at least two coordinates, they are forced to agree in all three. Thus, no two distinct elements agree in exactly two coordinates. 
\end{Proof}

\begin{Corollary}
    The independence number of the Cayley graph $Cay_{p^2}(G,C)$ is $\alpha^2 \beta^2 \gamma$. 
\end{Corollary}

\begin{Proof}
    We trivially construct the independent set $S$ as 
        \begin{align*}
            S = \bigcup_{(i,j,k) \in \mathcal{L}} C(i,j,k)
        \end{align*}
    where $\mathcal{L}$ is the set of indices previously defined in Lemma \ref{index in max construction}. 
\end{Proof}

After completing the analysis of the fundamental graph parameters, specifically the clique, chromatic, and independence numbers, we transition from examining vertex subsets to investigating the existence of spanning paths and cycles within the graph. The Hamiltonicity of a Cayley graph is a classical problem in algebraic graph theory. In the following proposition, we characterize the conditions under which $Cay_{p^2}(G,C)$ contains a Hamiltonian cycle or a Hamiltonian path, revealing a dependency on the smallest prime factor $\alpha$.

\begin{Proposition}
    The Cayley graph $Cay_{p^2}(G,C)$ is Hamiltonian for $\alpha =2$ and Semi-Hamiltonian for $\alpha>2$. 
\end{Proposition}

\begin{Proof}
    The vertex of $Cay_{p^2}(G,C)$ can be seen as the triple pair of the $\alpha$-, $\beta$-, and $\gamma$-components as follows.
        \begin{align*}
            V(Cay_{p^2}(G,C)) = \{ (g_\alpha, g_\beta, g_\gamma) : g_\alpha \in \mathbb{Z}_{\alpha^2}, g_\beta \in \mathbb{Z}_{\beta^2}, g_\gamma \in \mathbb{Z}_{\gamma^2} \}. 
        \end{align*}
    Two distinct vertices are adjacent if they share the same two coordinates and the third lie in a different residue class modulo $\alpha, \beta,$ or $\gamma$. So, for a fixed $g_\alpha \in \mathbb{Z}_{\alpha^2}$, all the vertices with the same $\alpha$-component can be listed as follow.
        \begin{align*}
            \begin{bmatrix}
                (\bar{0},\bar{0}) & (\bar{0},\bar{1}) & \cdots & (\bar{0}, \overline{\gamma^2 -1}) \\
                (\bar{1},\bar{0}) & (\bar{1},\bar{1}) & \cdots &(\bar{1}, \overline{\gamma^2 -1}) \\
                \vdots & \vdots & \ddots & \vdots \\
                (\overline{\beta^2 -1}, \bar{0}) & (\overline{\beta^2 -1}, \bar{1}) & \cdots & (\overline{\beta^2 -1}, \overline{\gamma^2 -1})
            \end{bmatrix}
        \end{align*}
    Now we construct a path traversing each vertex listed above. Start from the vertex $(\bar{0}, \bar{0})$ then go to the right traversing each vertex in the same row until the path reaches $(\bar{0}, \overline{\gamma^2 -1})$. Then, from $(\bar{0}, \overline{\gamma^2 -1})$, go down one step into $(\bar{1}, \overline{\gamma^2 -1})$. Continue traversing to the left, visiting each vertex in the row until it reaches $(\bar{1}, \bar{0})$. From $(\bar{1}, \bar{0})$ go down one step to the vertex $(\bar{2}, \bar{0})$. From there, repeat the step that is being used from $(\bar{0}, \bar{0})$, traversing to the right, and so on. Since there are odd number of rows, using the construction, the path should ends up in the vertex $(\overline{\beta^2 -1}, \overline{\gamma^2-1})$, visiting each vertex exactly once.

    Now, for $\alpha=2$, we have $g_\alpha \in \mathbb{Z}_4$, we can start our path from $g_\alpha \equiv \bar{0}$. Starting from the vertex $(\bar{0}, \bar{0}, \bar{0})$, using the construction before, we traverse all the vertex until it reaches $(\bar{0}, \overline{\beta^2 -1}, \overline{\gamma^2-1})$. From there, we then go to $(\bar{1},\overline{\beta^2 -1}, \overline{\gamma^2-1})$. Reversing the path construction as before, we reach $(\bar{1}, \bar{0}, \bar{0})$. We could do this procedure again and again until we eventually reach $(\bar{3}, \bar{0}, \bar{0})$, already visiting all of the vertices in the graph. Since $(\bar{3}, \bar{0}, \bar{0})$ is adjacent to $(\bar{0}, \bar{0}, \bar{0})$, we complete the path to go back to the original starting vertex. Thus, there exists a closed path that traverses each vertex exactly once, and therefore we conclude that $Cay_{p^2}(G,C)$ is a Hamiltonian graph.

    In the case that $\alpha>2$, using the same path construction, since there are odd number of elements in $\mathbb{Z}_{\alpha^2}$, then a path beginning from $(\bar{0}, \bar{0}, \bar{0})$ would end up in $(\overline{\alpha^2 -1}, \overline{\beta^2 -1}, \overline{\gamma^2 -1})$, after traversing all vertices in the graph. Because $(\bar{0}, \bar{0}, \bar{0})$ and  $(\overline{\alpha^2 -1}, \overline{\beta^2 -1}, \overline{\gamma^2 -1})$ are not adjacent, then the constructed path is not a closed path. And notice that since there are odd number of possible values in $\mathbb{Z}_{\alpha^2}$, any path starting from $(\bar{0}, \bar{0}, \bar{0})$ that visit each vertex exactly once could not end in a vertex that is adjacent to it. Hence, since we can already construct a path that traverses all of the vertex exactly once, we conclude that $Cay_{p^2}(G,C)$ is Semi-Hamiltonian. 
\end{Proof}

We provide the brief sketch of the Hamiltonian path for $\alpha =2, \beta=3$, and $\gamma=5$ for clearer understanding.

\begin{Example}
    We visualize all of the vertex in $Cay_{p^2}(G,C)$ where $G \cong \mathbb{Z}_{900}$ in a grid, similar to a coordinate in a cartesian $XYZ$-plane. We omit all the unecessary edges and only show the edges corresponding to the hamiltonian path.
    \begin{figure}[H]
        \centering
        \includegraphics[width=0.5\linewidth]{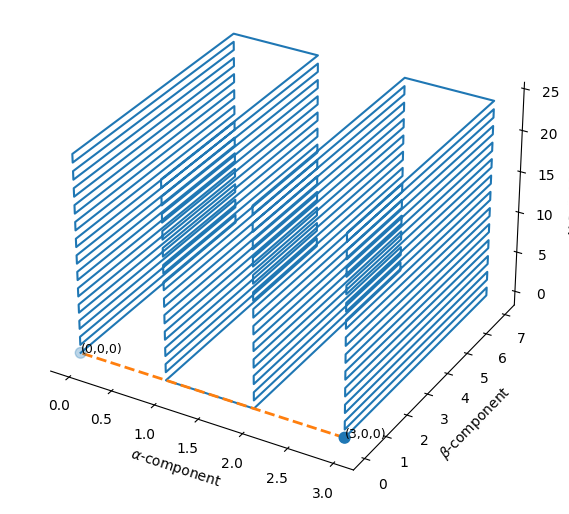}
        \caption{Hamiltonian path in $Cay_{p^2}(\mathbb{Z}_{900}, C)$}
        \label{fig:placeholder}
    \end{figure}
\end{Example}

Finally, we conclude our structural analysis by examining the distance properties of the graph. While previous propositions have described the internal density and the existence of spanning paths, the diameter provides a global measure of the graph's efficiency in terms of the maximum distance between any two vertices. 

\begin{Proposition}
    The diameter of the Cayley graph $Cay_{p^2}(G,C)$ is $6$. 
\end{Proposition}

\begin{Proof}
    Fix a $z \in \mathbb{Z}_{\gamma^2}$, then the subset
        \begin{align*}
            A = \{ (g_\alpha, g_\beta, z): g_\alpha \in \mathbb{Z}_{\alpha^2}, g_\beta \in \mathbb{Z}_{\beta^2} \}
        \end{align*}
    induces the Cayley graph $Cay(\mathbb{Z}_{\alpha^2 \beta^2}, S_A)$ where $S_A = \{ x \in G : |x| = \alpha^2 \text{ or } \beta^2 \}$. Moreover, it is clear that the distance between two vertices $(a,b), (x,y) \in V(Cay(\mathbb{Z}_{\alpha^2 \beta^2}, S_A))$ is as follows. For brevity, $ p \equiv q \text{ (mod }r)$ is denoted by $ p \equiv_r q$.
        \begin{align*}
            d_{Cay(\mathbb{Z}_{\alpha^2 \beta^2}, S_A)}((a,b),(x,y)) = \begin{cases}
                1, & a =x, b \equiv_{\beta} y \text{ or } b=y, a \equiv_\alpha x, \\
                2, & a \not\equiv_{\alpha} x \text{ and } b \not\equiv_{\beta} y, \\
                3, & a \equiv_{\alpha} x \text{ and } b \not\equiv_{\beta} y \text{ (or vice versa)} \\
                4, &a \equiv_{\alpha} x, \text{ and } b \equiv_{\beta} y.
            \end{cases}
        \end{align*}
    Now, consider two distinct vertices $(a,b,c), (x,y,z) \in V(Cay_{p^2}(G,C))$. Then
        \begin{enumerate}
            \item [$(i)$] If $c = z$, then $d_{Cay_{p^2}(G,C)}((a,b,c),(x,y,z)) \leq 4$. 
            \item [$(ii)$] If $c \neq z$, then either $c \equiv_{\gamma}z$ or $c \not\equiv_{\gamma} z$. Now, provided that $c \equiv_{\gamma} z$, we have
                \begin{align*}
                    d_{Cay_{p^2}(G,C)}((a,b,c),(x,y,z)) = \begin{cases}
                3, & a =x, b \equiv_{\beta} y \text{ or } b=y, a \equiv_\alpha x, \\
                4, & a \not\equiv_{\alpha} x \text{ and } b \not\equiv_{\beta} y, \\
                5, & a \equiv_{\alpha} x \text{ and } b \not\equiv_{\beta} y \text{ (or vice versa)} \\
                6, &a \equiv_{\alpha} x, \text{ and } b \equiv_{\beta} y.
            \end{cases}
                \end{align*}
            In the case that $c \not\equiv_{\gamma} z$, we obtain
                \begin{align*}
                    d_{Cay_{p^2}(G,C)}((a,b,c),(x,y,z)) = \begin{cases}
                2, & a =x, b \equiv_{\beta} y \text{ or } b=y, a \equiv_\alpha x, \\
                3, & a \not\equiv_{\alpha} x \text{ and } b \not\equiv_{\beta} y, \\
                4, & a \equiv_{\alpha} x \text{ and } b \not\equiv_{\beta} y \text{ (or vice versa)} \\
                5, &a \equiv_{\alpha} x, \text{ and } b \equiv_{\beta} y.
            \end{cases}
                \end{align*}
        \end{enumerate}
        Therefore, the diameter of the Cayley graph $Cay_{p^2}(G,C)$ is $6$.
\end{Proof}

\section{Conclusion}

In this paper, we have investigated the structural properties and graph parameters of the prime square order Cayley graph $Cay_{p^2}(G,C)$ for a cyclic group $G$ of order $\alpha^2 \beta^2 \gamma^2$, where $\alpha < \beta < \gamma$ are distinct prime numbers. By utilizing the cyclic decomposition of the group into $\mathbb{Z}_{\alpha^2} \times \mathbb{Z}_{\beta^2} \times \mathbb{Z}_{\gamma^2}$ , we have established the following key results on the regularity number, connectivity, eulerianity, clique number, independence number, planarity, hamiltonicity, and diameter as well.
These findings provide a comprehensive structural characterization of Cayley graphs for cyclic group of order $\alpha^2 \beta^2 \gamma^2$. For future research, this analytical approach could be extended to explore Cayley graphs of higher-order prime powers or investigated within the context of non-abelian groups.


\begin{thebibliography}{999}

\bibitem{cayley} { A. Cayley, The theory of groups: Graphical representation}, {\it Am. J. Math.}, \textbf{1 }(2) (1878), 174-176.

\bibitem{Graph theoretical concepts 2}
 S. Pirzada, {\em An introduction to graph theory}, Hyderabad, India: Universities Press Orient Blackswan,  2012.

\bibitem{Group theoretical concepts}
D. J. Robinson, {\em A Course in the Theory of Groups},
    New York-Heidelberg Berlin: Springer-Verlag:  1982.

\bibitem{shojaee} { I. Shojaee, A. Erfanian, B. Tolue}, { Some new approch on prime and composite order Cayley graphs}, {\it Quasigroups and Related Systems}, \textbf{27},($1$), (2019) $147-156$.

\bibitem{supar} T. Suparwati, Y, Susanti, S. Wahyuni, A. Erfanian, Prime Cubic Order Cayley Graph of Cyclic Groups, Asian-European Journal of Mathematics, \href{https://doi.org/10.1142/S1793557126500038}{https://doi.org/10.1142/S1793557126500038}    

\bibitem{YS} Y. Susanti and A. Erfanian,   {Prime Square Order Cayley Graph of Cyclic Groups}, {\it Asian European Journal of Mathematics} {\bf 17} (2), (2024) \href{https://doi.org/10.1142/S1793557124500037}{https://doi.org/10.1142/S1793557124500037}. 

\bibitem{tolue1} {B. Tolue}, { Some graph parameters on the composite order Cayley graph}, {\it Caspian J. Mathematical Sciences}, \textbf{8} (1), (2019) $10-17$.

\bibitem{tolue2} { B. Tolue}, { The prime order Cayley graph}, {\it U. P. B Sci. Bull., Series A.} \textbf{77},($3$), (2015) $207-218$.


\bibitem{west}
     D. B. West, {\em Introduction to Graph Theory},      2 Eds., USA: Prentice Hall,  2001.
     






\end{thebibliography}
\end{document}